\documentclass[10pt]{article}
\usepackage{}
\usepackage{amsfonts}
\usepackage{mathrsfs}
\usepackage{amssymb}
\usepackage{amsmath,amsthm,amscd,amssymb,mathrsfs,graphicx,amsfonts,mathrsfs}
\usepackage{amsmath}
\usepackage[all]{xy}

\newtheorem{theorem}{\hspace{0.7cm}Theorem }[section]
\newtheorem{corollary}[theorem]{\hspace{0.7cm}Corollary }
\newtheorem{lemma}[theorem]{\hspace{0.7cm}Lemma}
\newtheorem{proposition}[theorem]{\hspace{0.7cm}Proposition }

\newtheorem{remark}[theorem]{\hspace{0.7cm}Remark }

\def\mod{\mathop{\rm mod}\nolimits}

\def\Hom{\mathop{\rm Hom}\nolimits}

\def\id{\mathop{\rm id}\nolimits}

\def\Gproj{\mathop{\rm Gproj}\nolimits}
\def\Ginj{\mathop{\rm Ginj}\nolimits}
\def\proj{\mathop{\rm proj}\nolimits}
\def\inj{\mathop{\rm inj}\nolimits}

\def\Id{\mathop{\rm Id}\nolimits}

\title{\large \bf Auslander-Reiten Triangles in Homotopy Categories
}
\author{Yuefei Zheng\thanks{{\it E-mail address}: yuefeizheng@sina.com}
\ and Zhaoyong Huang\thanks{{\it E-mail address}: huangzy@nju.edu.cn}\\
{\footnotesize \it Department of Mathematics, Nanjing
University, Nanjing 210093, Jiangsu Province, China}}
\date{}
\begin{document}
\baselineskip=16pt \maketitle

\begin{abstract}
Let $A$ be an artin algebra. We show that the bounded homotopy category of finitely generated right $A$-modules
has Auslander-Reiten triangles. Two applications are given: (1) we provide an alternative proof of a theorem of Happel in [H2];
(2) we prove that over a Gorenstein algebra, the bounded homotopy category of finitely generated Gorenstein projective
(resp. injective) modules admits Auslander-Reiten triangles, which improves a main result in [G].
\end{abstract}

\bigskip
\section{ Introduction}
%\lbl{s1.1}
\setcounter{equation}{0}

Throughout this paper, $A$ is an artin algebra over a fixed commutative artin ring $R$ and $D:=\Hom_{R}(-,E(R/J))$ is the usual duality,
where $J$ is the Jacobson radical of $R$ and $E(R/J)$ is the injective envelope of $R/J$. We denote by $\mod A$ the category of
finitely generated right $A$-modules. As usual, we write $\proj A$ (resp. $\inj A$) the category of finitely generated projective
(resp. injective) right $A$-modules.

Auslander-Reiten sequences, also known as almost split sequences, are one of the central tools in the representation theory of
artin algebras. Auslander-Reiten triangles can also be defined by almost split morphisms in $\Hom$-finite Krull-Schmidt triangulated $R$-categories.
Happel proved in [H3] that the bounded homotopy category $K^b(\proj A)$ of finitely generated projective right $A$-modules has
right Auslander-Reiten triangles if and only if the left self-injective
dimension of $A$ is finite; and dually, the bounded homotopy category $K^b(\inj A)$ of finitely generated injective
right $A$-modules has left Auslander-Reiten triangles if and only if the right self-injective dimension of $A$ is finite.

There is a close relation between Auslander-Reiten triangles and Serre functors ([RV]):
a $\Hom$-finite Krull-Schmidt triangulated $R$-category has right (resp. left)
Auslander-Reiten triangles if and only if it has a right (resp. left) Serre functor.
A Serre functor by definition is a right Serre functor which is an equivalence.
In [BJ], Backelin and Jaramillo proved that the bounded homotopy category
$K^b(\mod A)$ of $\mod A$ has a right Serre duality.
Their method is based on the construction of a $t$-structure in $K^b(\mod A)$,
and their proof is somewhat complicated although they obtained some more results. We use the terminology of Auslander-Reiten
triangles to prove that the right Serre functor in $K^b(\mod A)$ is always an equivalence (Theorem 3.4).
Our result is based on the fact that right (left) minimal almost split morphisms are stable under quotients.
It seems more elementary. In particular, we determine Auslander-Reiten triangles admitting special ending (resp. starting) terms (Proposition 4.2).

As in abelian categories, the existence of Auslander-Reiten triangles in subcategories were investigated by J\o rgensen in [J].
Note that $K^{b}(\proj A)$ and $K^{b}(\inj A)$ can be embedded in $K^{b}(\mod A)$ naturally. By using the obtained result
about the existence of Auslander-Reiten triangles in $K^b(\mod A)$, we reprove the Happel's theorem mentioned above (Theorem 4.4).
The steps we take seem more ``categorization" and can be easily treated dually. The advantage here is that the Auslander-Reiten triangles
we treat always lie in $K^{b}(\mod A)$. Similarly, we prove that over a Gorenstein algebra, the bounded homotopy category of
finitely generated Gorenstein projective (resp. injective) modules admits Auslander-Reiten triangles, which improves a main result in [G].

\section{Preliminaries}
%\lbl{s1.1}
\setcounter{equation}{0}

\vspace{0.2cm}

Recall that a right $A$-module $M$ is called {\it Gorenstein projective} if there exists an exact sequence $T^{\bullet}$
of projective modules which remains exact when applying the functor $\Hom_A(-,P)$ for any $P\in \proj A$ such that $M$
is isomorphic to some kenrel of $T^{\bullet}$. Dually, The notion of {\it Gorenstein injective modules} is defined.
We denote the category of all finitely generated Gorenstein projective (resp. injective) modules by $\Gproj A$
(resp. $\Ginj A$). Note that we have $\proj A \subset \Gproj A$ and $\inj A \subset \Ginj A$.

Let $f:M\to N$ be a morphism in $\mod A$. According to [AR1], $f$ is called {\it right almost split} if it is not
a retraction, and any morphism $g:L\to N$ which is not a retraction factors through $f$; it is called {\it right minimal} if
any morphism $h$ satisfying $f=f\cdot h$ is an automorphism of $M$; and it is called {\it right minimal almost split} if it is both
right almost split and right minimal. The left versions are defined dually. An exact sequence
$$0\to M \overset{\alpha} \to N \overset{\beta} \to L\to 0$$
in $\mod A$ is called an {\it almost split sequence} if either $\alpha$ is left minimal almost split or $\beta$ is right minimal almost split.
This also means that $M$ and $L$ must be indecomposable. We say $\mod A$ has {\it almost split sequences} if for any indecomposable
non-injective module $M$ in $\mod A$ there is an almost split sequence starting at $M$,
and for any indecomposable non-projective module $N \in \mod A$
there is an almost split sequence ending at $N$. By [AR1], $\mod A$ always has almost split sequences.

Let $\mathcal{T}$ be
a $\Hom$-finite Krull-Schmidt triangulated $R$-category. The notion of almost split triangles in $\mathcal{T}$ was introduced by
Happel in [H1]. A triangle
$$X \overset{\alpha}\to Y \overset{\beta}\to Z \to X[1]$$
in $\mathcal{T}$ is called an {\it Auslander-Reiten triangle} (or {\it almost split triangle}) if either $\alpha$
is left minimal almost split or $\beta$ is right minimal almost split (see [H2], where the notions source morphisms and
sink morphisms were used). A triangulated category is said to {\it have right $($resp. left$)$ Auslander-Reiten triangles}
if for any indecomposable object $M$ there is an Auslander-Reiten triangle ending (resp. starting) at $M$.

Let $F:\mathcal{T} \to \mathcal{T}$ be a triangulated functor. According to [RV], $F$ is called a {\it right Serre functor} if
for any $X$, $Y$ in $\mathcal{T}$, there is an isomorphism $D\Hom_{\mathcal{T}}(X,Y)\cong \Hom_{\mathcal{T}}(Y,FX)$ which
is natural in $X$ and $Y$. This $F$ is unique up to natural isomorphism. A left Serre functor is defined dually.

\begin{theorem}$\mathrm{([BJ, Corollary~2.5 ~and~ Proposition~ 4.6])}$

Let $A$ be an Artin algebra. Then there is a right serre functor $S: K^{b}(\mod A)\to K^{b}(\mod A)$,
equivalently, $K^{b}(\mod A)$ has right Auslander-Reiten triangles.

\end{theorem}

\bigskip
\section{ AR-triangles in homotopy categories}
%\lbl{s1.1}
\setcounter{equation}{0}

Let $\mathcal{A}$ be an additive $R$-category. By an $\mathcal{A}$-module, we mean a contravariant $R$-linear
functor from $\mathcal{A}$ to the category of $R$-modules. We denote by $\mod \mathcal{A}$ the category of
finitely presented $\mathcal{A}$-modules. Note that, in general, $\mod \mathcal{A}$ is not an abelian category.
It is an abelian category if and only if $\mathcal{A}$ has pseudo-kernels ([A]). We call $\mathcal{A}$ a {\it dualizing $R$-category}
if $D$ gives a duality between $\mod \mathcal{A}$ and $\mod \mathcal{A}^{op}$. Note that, in this case,
$\mod \mathcal{A}$ is always an abelian category, and hence the bounded complex category $C^{b}(\mod \mathcal{A})$
of $\mod \mathcal{A}$ is also an abelian category. We begin with a main theorem in [BJR].

\begin{lemma}$\mathrm{([BJR, Theorem~4.3])}$
Let $\mathcal{A}$ be a dualizing $R$-category. Then $C^{b}(\mod \mathcal{A})$ has almost split sequences.
\end{lemma}

Note that $\mod A$ is equivalent to $\mod (\proj A)$ as additive $R$-categories ([A]). This means that $C^{b}(\mod A)$ always
has almost split sequences. As a consequence, for any indecomposable non-projective (resp. non-injective) object
$X$ in $C^{b}(\mod A)$, there is always an almost split sequence ending (resp. starting) at $X$.

Let $\mathcal{B}$ be an additive category and $\mathcal{C}$ an additive full subcategory of $\mathcal{B}$ closed under summands.
Then we can form the factor category $\mathcal{B}/\mathcal{C}$. The objects in $\mathcal{B}/\mathcal{C}$
are the same as in $\mathcal{B}$, and the morphisms are the morphisms in $\mathcal{B}$ modulo morphisms factor
through an object of $\mathcal{C}$. There is a natural factor functor $\pi: \mathcal{B}\to \mathcal{B}/\mathcal{C}$.
It is an additive functor. For both objects and morphisms, we denote their images under $\pi$ by adding $\ \widetilde{ }\ $ above.
The following lemma is a well known fact. For the reader's convenience we give a quick proof here.

\begin{lemma}
Let $\mathcal{B}$ be a Hom-finite Krull-Schmidt $R$-category and $\mathcal{C}$ an additive full subcategory of $\mathcal{B}$
closed under summands. If $f: M\to N$ is a right (resp. left) minimal almost split morphism in $\mathcal{B}$, then
$\widetilde{f}: \widetilde{M}\to \widetilde{N}$ is also a right (resp. left) minimal almost split morphism in $\mathcal{B}/\mathcal{C}$.
\end{lemma}

{\it Proof.} Obviously, $N$ is indecomposable; in particular, it has no nonzero summands in $\mathcal{C}$. By [AR2, Lemma 1.1(c)],
$\widetilde{f}: \widetilde{M}\to \widetilde{N}$ is not a retraction. Let $\widetilde{g}: \widetilde{L}\to \widetilde{N}$
be not a retraction. Then it is induced by a morphism $g: L\to N$ which is also not a retraction. So $g$ factors through $f$
since $f$ is right almost split, and hence $\widetilde{g}$ factors through $\widetilde{f}$. If $\widetilde{f}$ is not minimal,
then $\widetilde{f}$ has a direct summand of the form $\widetilde{W} \to 0$, and so $f$ has a direct summand of the
form $W \to 0$ or $W \to C$ with $0 \neq C\in \mathcal{C}$. Note that the former one gives a contradiction to the minimality of $f$,
and the latter one gives a contradiction to the indecomposableness of $N$. \hfill$\square$

\vspace{0.2cm}

A complex $X$ is called {\it contractible} if it is isomorphic to zero in $K^{b}(\mod A)$, that is,
it is splitting exact. Note that a chain map of complexes is homotopic to zero if and only if it factors through some contractible
complex. So $K^{b}(\mod A)$ is exactly the factor category of $C^{b}(\mod A)$ by modulo contractible complexes ([H2, p.28]).
We also need the following

\begin{lemma}
A complex $X$ is a projective (resp. injective) object in $C^{b}(\mod A)$ if and only if it is a contractible
complex consisting of projective (resp. injective) modules in $\mod A$.
\end{lemma}
{\it Proof.} See for example [EJ2, Theorem 1.4.7]. \hfill$\square$

\vspace{0.2cm}

Now we can prove the following result, which improves Theorem 2.1.

\begin{theorem}
$K^{b}(\mod A)$ has Auslander-Reiten triangles.
\end{theorem}

{\it Proof.} Let $0\neq \widetilde{X} \in K^{b}(\mod A)$ be indecomposable. Then it is induced by an indecomposable object
$X$ in $C^{b}(\mod A)$. Note that $X$ is neither projective nor injective by Lemma 3.3. It follows from Lemma 3.1 that
there is a minimal right almost split morphism $f:Y\to X$ in $C^{b}(\mod A)$. Then by Lemma 3.2, its image $\widetilde{f}:
\widetilde{Y}\to \widetilde{X}$ is also right minimal almost split. Complete it to a triangle
$$L \to\widetilde{Y}\to \widetilde{X} \to L[1]$$
in $K^{b}(\mod A)$. Then by definition, it is an Auslander-Reiten triangle in $K^{b}(\mod A)$. Dually,
there is an Auslander-Reiten triangle starting at $\widetilde{X}$. \hfill$\square$

\section{Applications}

In this section, we will reprove the Happel's theorem by using the main result in Section 3. Our proof is
based on the restriction of Auslander-Reiten triangles in subcategories. Also one can see that using this
technique, over a Gorenstein algebra the existence of Auslander-Reiten triangles in the bounded homotopy
category of Gorenstein projective modules is valid. This improves a result by Gao in [G] where only the
existence of right Auslander-Reiten triangles is proved and the condition CM-finiteness is necessary.

Although $K^{b}(\mod A)$ has Auslander-Reiten triangles, it is difficult to compute the Auslander-Reiten translation.
By [RV, Theorem I.2.4], $K^{b}(\mod A)$ has Auslander-Reiten triangles if and only if it has a Serre functor.
Denote the Serre functor of $K^{b}(\mod A)$ by $S$. Then by [RV, Proposition I.2.3], we have that
for any indecomposable object $X$,
the left end term of an Auslander-Reiten triangle ending at $X$ is $S\cdot [-1]$, the other end term is its quasi-inverse
$S^{-1}\cdot [1]$. Thanks to this result, we only need to compute the Serre dual object for an indecomposable object $X$.

\begin{lemma}Let $X$ and $Y$ be in $C^{b}(\mod A)$.
\begin{enumerate}
\item[(1)] If $X$ is degreewise projective, then we have a natural isomorphism
$$D\Hom_{A}(X,Y)\cong \Hom_{A}(Y,X\otimes_{A}DA).$$

\item[(2)] If $Y$ is degreewise injective, then we have a natural isomorphism
$$D\Hom_{A}(X,Y)\cong \Hom_{A}(\Hom_{A}(DA,Y),X).$$
\end{enumerate}
\end{lemma}

{\it Proof.} (1) Note that for any $X,Y$ in $\mod A$ , we have a natural morphism
$$\delta_{Y,X}:Y\otimes_A\Hom_A(X,A)\to \Hom_A(X,Y),~y\otimes f\longmapsto (x\longmapsto yf(x)).$$
So we have a natural morphism
$$\eta_{X,Y}:D\Hom_A(X,Y)\overset{D\delta_{Y,X}}\longrightarrow D(Y\otimes_A\Hom_A(X,A))$$
\begin{align*}
&~~~~~~~~~~~~~~~~~~~~~~~~~~~~~~~~~~~~~~~~~~~~~~~~~~~~~~~~\cong\Hom_R(Y\otimes_A\Hom_A(X,A),E(R/J))\\
&~~~~~~~~~~~~~~~~~~~~~~~~~~~~~~~~~~~~~~~~~~~~~~~~~~~~~~~~\cong \Hom_A(Y,\Hom_R(\Hom_A(X,A),E(R/J))\\
&~~~~~~~~~~~~~~~~~~~~~~~~~~~~~~~~~~~~~~~~~~~~~~~~~~~~~~~~\cong \Hom_A(Y,X\otimes_ADA).
\end{align*}
It is known that if $X$ is projective, then $\delta_{Y,X}$ is an isomorphism. Thus $D\delta_{Y,X}$ is also an isomorphism.
Therefore we have a natural isomorphism
$$\eta_{X,Y}:D\Hom_{A}(X,Y)\cong \Hom_{A}(Y,X\otimes_{A}DA).$$Now the isomorphism can be extended to the desired situation.

(2) Let $Y$ be injective. Then $Y\cong DA\otimes_A P$ for some $P\in \proj A$. Then we have isomorphisms
\begin{align*}
&\ \ \ \ \Hom_{A}(\Hom_{A}(DA,Y),X)\\
& \cong \Hom_{A}(\Hom_{A}(DA,DA\otimes_A P),X)\\
& \cong \Hom_A(P,X)\\
& \cong D\Hom_{A}(X,P\otimes_{A}DA)\ {\rm (by\ (1))}\\
& \cong D\Hom_{A}(X,Y).
\end{align*}
Similarly, the isomorphism can be extended to the desired situation.
\hfill$\square$

\begin{proposition}Let $X$ and $Y$ be in $C^{b}(\mod A)$.
\begin{enumerate}
\item[(1)] If $X$ is degreewise projective, then there is an Auslander-Reiten triangle in $K^{b}(\mod A)$
$$X[-1]\otimes _A DA\to M\to X\to X\otimes _A DA.$$

\item[(2)] If $Y$ is degreewise injective, then there is an Auslander-Reiten triangle in $K^{b}(\mod A)$
$$Y\to N \to \Hom_A(DA,Y[1])\to Y[1].$$
\end{enumerate}
\end{proposition}

{\it Proof.} We only prove (1), and the proof of (2) is similar. Let $X \in K^{b}(\proj A)$. By Lemma 4.1, we have
$D\Hom_{A}(X,Y)\cong \Hom_{A}(Y,X\otimes_{A}DA)$. Thus we have isomorphisms
$$D\Hom_{K^{b}(\mod A)}(X,Y)\cong DH^{0}\Hom_{A}(X,Y) \cong H^{0}D\Hom_{A}(X,Y)\cong $$
$$ H^{0}\Hom_{A}(Y,X\otimes_{A}DA) \cong \Hom _{K^{b}(\mod A)}(Y,X\otimes_{A}DA).$$ This holds for any $Y\in K^{b}(\mod A)$.
Then by the Yoneda's lemma, the Serre dual object for $X$ is $X\otimes_{A}DA$. Now by [RV, Proposition I.2.3],
we have the desired triangle. \hfill$\square$

\vspace{0.2cm}

Let $\mathcal{B}$ be an additive category and $\mathcal{C}$ a full subcategory of $\mathcal{B}$.
Recall that a morphism $f:B\to C$ with $B\in \mathcal{B}$ and $C\in \mathcal{C}$ is called a {\it $\mathcal{C}$-preenvelope}
if the natural map $\Hom_{\mathcal{B}}(C,C')\to \Hom_{\mathcal{B}}(B,C')\to 0$ is exact for any
$C'\in \mathcal{C}$. A $\mathcal{C}$-preenvelope $C$
is called a {\it $\mathcal{C}$-envelope} if it is left minimal. Dually, the notion of (pre)covers is defined ([AR3, E]).
The following proposition involves Auslander-Reiten triangles in subcategories.

\begin{lemma}$\mathrm{([J, Theorem~3.1~\mathrm{and}~Theorem~3.2])}$ Let $\mathcal{T}$ be a triangulated category
and $\mathcal{C}$ a full subcategory of $\mathcal{T}$ closed under extensions.
\begin{enumerate}
\item[(1)] Let $X \to Y \to Z \to X[1]$ be an Auslander-Reiten triangle in $\mathcal{T}$ with $Z \in \mathcal{C}$.
If there is an object $A' \in \mathcal{C} $ with a nonzero morphism $Z\to A'[1]$, then the following are equivalent.

~~~~$\bullet$ $X$ has a $\mathcal{C}$-cover of the form $A\to X$.

~~~~$\bullet$ There is an Auslander-Reiten triangle $A \to B \to Z \to A[1]$ in $\mathcal{C}$.

\item[(2)] Let $X \to Y \to Z \to X[1]$ be an Auslander-Reiten triangle in $\mathcal{T}$ with $X \in \mathcal{C}$.
If there is an object $Z' \in \mathcal{C} $ with a nonzero morphism $Z'\to X[1]$, then the following are equivalent.

~~~~$\bullet$ $Z$ has a $\mathcal{C}$-envelope of the form $Z\to N$.

~~~~$\bullet$ There is an Auslander-Reiten triangle $X \to M \to N \to X[1]$ in $\mathcal{C}$.
\end{enumerate}
\end{lemma}

As an application of Theorem 3.4, we now are in a position to reprove the following Happel's theorem.
Our argument is very different from the original one.

\begin{theorem}$\mathrm{([H2, Section~3.4])}$
\begin{enumerate}
\item[(1)] $K^{b}(\proj A)$ has right Auslander-Reiten triangles if and only if $\id A_{A^{op}}<\infty$.
\item[(2)] $K^{b}(\inj A)$ has left Auslander-Reiten triangles if and only if $\id A_{A}<\infty$.
\end{enumerate}
\end{theorem}

{\it Proof.} We only prove (2), and (1) is its dual.

Let $0\neq Y\in K^{b}(\inj A)$ be indecomposable. Then by Proposition 4.2,
there is an Auslander-Reiten triangle
$$Y\to L \to \Hom_{A}(DA,Y)[1] \to Y[1]$$
in $K^{b}(\mod A)$. Since $(Y[1])[-1]\overset{\Id_{Y}}\to Y$ is not homotopic to zero, by Lemma 4.3
we have that $K^{b}(\inj A)$ has left Auslander-Reiten triangles if and only if $\Hom_{A}(DA,Y)$ has a
$K^{b}(\inj A)$-envelope for any $Y \in K^{b}(\inj A)$. Now it suffices to prove that $\id A_{A}<\infty$
if and only if $\Hom_{A}(DA,Y)$ has a $K^{b}(\inj A)$-envelope for any $Y \in K^{b}(\inj A)$.

If $\id A_{A}<\infty$, then the injective dimension of any module in $\proj A$ is finite.
Since $\Hom_{A}(DA,Y)$ consists of modules in $\proj A$, then by using induction on
the width of $\Hom_{A}(DA,Y)$, one can get an injective resolution $f:\Hom_{A}(DA,Y)\to L$, where $f$ is
a quasi-isomorphism and $L\in K^{b}(\inj A)$. If there is a chain map $\alpha :\Hom_{A}(DA,Y)\to I$
with $I\in K^{b}(\inj A)$, then $ f^{\ast}:\Hom_{A}(L,I) \to \Hom_{A}(\Hom_{A}(DA,Y),I)$ is a quasi-isomorphism
and $H^{0}(f^{\ast})$ is an isomorphism. Hence there is some $\beta:L\to I$
such that $\alpha$ homotopic to $\beta \cdot f$.
If there is some $g$ satisfying $g\cdot f$ homotopic to $f$, then it is a quasi-isomorphism,
and hence a homotopic equivalence. It follows that $L$ is a $ K^{b}(\inj A)$-envelope of $\Hom_{A}(DA,Y)$.

Now suppose that $\Hom_{A}(DA,Y)$ has a $K^{b}(\inj A)$-envelope for any $Y \in K^{b}(\inj A)$.
Let $Y=DA$. Then $\Hom_{A}(DA,Y)\cong A$. Let $A\to I$ be the $K^{b}(\inj A)$-envelope of $A$.
Complete it to a triangle
$$A\overset{\alpha}\to I \to L \to A[1]$$
in $K^{b}(\mod A)$. Then we have that $\Hom_{K^{b}(\mod A)}(L,Z)\cong 0$ for any $Z\in K^{b}(\inj A)$
by the Wakamatsu lemma (see for example [J, Lemma 2.1]); in particular, $\Hom_{K^{b}(\mod A)}(L,(DA)[i])=0$ for all $i$,
which implies that $L$ is exact. As a consequence, $\alpha$ is a quasi-isomorphism.
Since any injective resolution of $A$ is homotopically equivalent to $I$.
It follows that $\id A_{A}<\infty$. \hfill$\square$

\begin{remark}
Note that the functor $-\otimes_{A}DA: K^{b}(\proj A)\to K^{b}(\inj A)$ is
an equivalence. Hence the Gorenstein symmetry conjecture, which states that $\id A _{A}=\id A_{A^{op}}$ for any artin algebra $A$,
can be reformulated as follows.

$\bullet K^{b}(\proj A)$ has right Auslander-Reiten triangles if and only if it has left Auslander-Reiten triangles.
Dually,

$\bullet K^{b}(\inj A)$ has right Auslander-Reiten triangles if and only if it has left Auslander-Reiten triangles.
\end{remark}

\begin{corollary} The following are equivalent.
\begin{enumerate}
\item[(1)] $A$ is a Gorenstein algebra.
\item[(2)] $K^{b}(\proj A)$ has Auslander-Reiten triangles.
\item[(3)]$K^{b}(\inj A)$ has Auslander-Reiten triangles.
\end{enumerate}
\end{corollary}

Let $A$ and $B$ be artin algebras. According to [R], $A$ and $B$ are derived equivalent if and only if
$K^{b}(\proj A)$ and $K^{b}(\proj B)$ are equivalent as triangulated categories.

\begin{corollary}
Let $A$ and $B$ be artin algebras. If $A$ and $B$ are derived equivalent,
then $A$ is Gorenstein if and only if $B$ is Gorenstein.
\end{corollary}

In general, Theorem $3.4$ only tells us the validity of Auslander-Reiten triangles in $K^{b}(\mod A)$.
When consider some subcategory $\mathcal{C}$ of $K^{b}(\mod A)$, one usually relies on the restriction as
in Lemma $4.3$. However, it is often difficult to compute the serre dual objects for $\mathcal{A}$.
For example, the isomorphism in Lemma 4.2 for projective modules can not be extended to Gorenstein
projective version unless $A$ is self-injective, see [ASS, Lemma 2.12]. In the following, we only
consider the subcategory $K^{b}(\Gproj A$) (resp. $K^{b}(\Ginj A$)).

It was proved in [EJ1] that over a Gorenstein algebra any finitely generated module admits
a finitely generated Gorenstein projective precover. That is, for any $M\in \mod$ A, there
is a complex $G_M$ consisting of modules in $\Gproj A$ and a chain map $G_M\to M$ which is
a quasi-isomorphism after applying the functor $\Hom_A(G',-)$ for any $G'\in \Gproj A$.
Since $M$ has finite Gorenstein projective dimension, $G_M$
can be selected to be in $K^{b}(\Gproj A)$  by [Ho, Proposition 2.18].
The dual version for finitely generated Gorenstein injective modules is also valid.

\begin{theorem}
Let $A$ be a Gorenstein algebra. Then $K^{b}(\Gproj A)$ has Auslander-Reiten triangles.
\end{theorem}

{\it Proof.} The proof is similar to that of Theorem 4.4. First, we prove that $K^{b}(\Gproj A)$
has right Auslander-Reiten triangles. Let $0\neq X\in K^{b}(\Gproj A)$ be indecomposable.
Then by Theorem 3.4, we have an Auslander-Reiten triangle
$$Y\to L\to X\to Y[1]$$
in $K^{b}(\mod A)$. We only need to prove that $Y$ has a $K^{b}(\Gproj A)$-cover.
In fact, we will prove that any $Y\in K^b(\mod A)$ has a $K^{b}(\Gproj A)$-cover. Note that
for any $M\in\mod A$, there is a chain map $G_M\to M$ with $G_M\in K^{b}(\Gproj A)$, which
is a quasi-isomorphism after applying the functor $\Hom_A(G',-)$ for any $G'\in \Gproj A$ as above.
By using induction on the width of $Y$, we have that there is a chain map $f_Y:G_Y\to Y$ with
$G_Y\in K^b(\Gproj A)$, which is also a quasi-isomorphism after applying the functor $\Hom_A(G',-)$
for any $G'\in \Gproj A$. Hence $\Hom_A(G',f_Y)$ is also a quasi-isomorphism for any
$G'\in K^{b}(\Gproj A)$. It is easy to see that $G_Y$ is a $K^{b}(\Gproj A)$-cover of $Y$.
If we consider the category $K^{b}(\Ginj A)$, we then obtain that $K^{b}(\Ginj A)$ admits
left Auslander-Reiten triangles. Note that $-\otimes_{A}DA:K^b(\Gproj A)\to K^b(\Ginj A)$
is an equivalence by [B]. This completes the proof.\hfill$\square$

\vspace{0.5cm}

\noindent $\mathbf{Acknowledgements}$

This research was partially supported by NSFC (Grant Nos. 11171142 and 11571164)
and an Innovation Programme of Jiangsu province (Grant No. 0203001505). The authors would like to express
their sincere thanks to Professor Silvana Bazzoni and the referees
for many considerable suggestions, which have greatly improved this paper.


\begin{thebibliography}{101}

\bibitem[ASS]{A1} I. Assem, D. Simson and A. Skowro$\acute{n}$ski, Elements of the representation theory of associative algebras.
Vol. 1. London Mathematical Society Student Texts {\bf 65}, Cambridge
University Press, Cambridge, 2006.

\bibitem[A]{A1} M. Auslander, {\it Representation Dimension of Artin Algebras}, Queen Mary College Mathematics Notes, 1971.
Republished in: Selected Works of Maurice Auslander, Amer. Math. Soc., Providence, 1999.

\bibitem[AR1]{A1} M. Auslander and I. Reiten, {\it Representation theory of Artin Algebras III: Almost split sequences}, Comm. Algebra {\bf 3} (1975), 239--294.

\bibitem[AR2]{A1} M. Auslander and I. Reiten, {\it Representation theory of Artin algebras V. Methods for computing almost split sequences and irreducible morphisms},
Comm. Algebra {\bf 5} (1977), 519--554.

\bibitem[AR3]{A1} M. Auslander and I. Reiten, {\it Applications of contravariantly finite subcategories}. Adv. Math. {\bf 86} (1991), 111--152.

\bibitem[BJ]{A1} E. Backelin and O. Jaramillo, {\it Auslander-Reiten sequences and t-structures on the homotopy category of an abelian category},
J. Algebra {\bf 339} (2011), 80--96.

\bibitem[BJR]{A1} R. Bautista, M. Jos\'e Souto Salorio and R. Zuazua, {\it Almost split sequences for complexes of fixed size}, J. Algebra {\bf 287} (2005), 140--168.

\bibitem[B]{A1} A. Beligiannis, {\it Cohen-Macaulay modules, (co)torsion pairs and virtually Gorenstein algebras}, J. Algebra {\bf 288} (2005), 137--211.

\bibitem[E]{A1} E.E. Enochs, {\it Injective and flat covers, envelopes and resolvents}, Israel J. Math. {\bf 39} (1981), 189--209.

\bibitem[EJ1]{A1} E.E. Enochs and O.M. Jenda, Relative Homological Algebra, Vol. 1.
Second revised and extended edition. de Gruyter Expositions in Mathematics {\bf 30},
Walter de Gruyter GmbH $\&$ Co. KG, Berlin, 2011.

\bibitem[EJ2]{A1} E.E. Enochs and O.M. Jenda, Relative Homological Algebra, Vol. 2. de Gruyter Expositions in Mathematics {\bf 54},
Walter de Gruyter GmbH $\&$ Co. KG, Berlin, 2011.

\bibitem[G]{A1} G. Nan, {\it Auslander-Reiten triangles on Gorenstein derived categories}, Comm. Algebra {\bf 40} (2012), 3912--3919.

\bibitem[H1]{A1} D. Happel, {\it On the derived category of a finite-dimensional algebra}, Comment. Math. Helvetici {\bf 62} (1987), 339--389.

\bibitem[H2]{A1} D. Happel, Triangulated categories in the representation theory of finite-dimensional algebras, London Mathematical Society Lecture Note Series, {\bf 119}, Cambridge University Press, Cambridge, 1988.

\bibitem[H3]{A1} D. Happel, {\it On Gorenstein algebras}, Representation theory of finite groups and finite-dimensional algebras (Bielefeld, 1991), 389--404,
Progr. Math., {\bf 95}, Birkh$\mathrm{\ddot{a}}$user, Basel, 1991.

\bibitem[Ho]{A1} H. Holm, {\it Gorenstein homological dimensions}, J. Pure Appl. Algebra {\bf 189} (2004),
167--193.

\bibitem[J]{A1} P. J\o rgensen, {\it Auslander-Reiten triangles in subcategories}, J. K-Theory {\bf 3} (2009), 583--601.

\bibitem[R]{A1} J. Rickard, {\it Morita theory for derived categories}, J. London Math. Soc. {\bf 39}  (1989), 436--456.

\bibitem[RV]{A1} I. Reiten and M. van den Bergh, {\it Noetherian hereditary abelian categories satisfying Serre duality}, J. Amer. Math. Soc. {\bf 15} (2002), 295--366.

\end{thebibliography}
\end{document}